\documentclass[a4paper,14pt]{article}
\usepackage[cp1251]{inputenc}
\usepackage[english,russian]{babel}
\usepackage{amsfonts,amssymb,mathrsfs,amscd}
\usepackage[tbtags]{amsmath}
\usepackage{amsfonts,amssymb,mathrsfs,amscd}
\usepackage{graphics}
\usepackage{wrapfig}
\usepackage{euscript}
\usepackage[dvips]{graphicx}
\usepackage{xcolor}

\date{}
\newtheorem{theorem}{Theorem}[section]
\newtheorem{lemma}{Lemma}[section]
\newtheorem{definition}{Definition}[section]

\begin{document}
\renewcommand{\baselinestretch}{0.96}
\renewcommand{\thesection}{\arabic{section}}
\renewcommand{\theequation}{\thesection.\arabic{equation}}
\csname @addtoreset\endcsname{equation}{section}
\Large

 \centerline{\bf On approximation by  tight wavelet frames on the field of $p$-adic numbers}
  \vskip0.5cm
  \centerline{\bf S.\,F.~Lukomskii\footnote[1]{Corresponding author, \\
  E-mail addresses: LukomskiiSF@info.sgu.ru (S.Lukomskii ),
 vam21@yandex.ru (A.Vodolazov)}, A.M.Vodolazov}
\vskip0.5cm
\noindent
 N.G.\ Chernyshevskii Saratov State University, 83 Astrakhanskaya Street, Saratov, 410012, Russian Federation

\begin{abstract} \large
We discuss the problem on approximation by  tight step wavelet frames on the field $\mathbb{Q}_p$ of $p$-adic numbers.    Let $G_n=\{x=\sum_{k=n}^\infty x_k p^k\}$, $X$ be a set of characters.  We define a step function $\lambda({\chi})$ that is constant on cosets  ${G}_n^\bot\setminus{G}_{n-1}^\bot$ by equalities
$\lambda ({G}_n^\bot\setminus{G}_{n-1}^\bot)=\lambda_n>0$
for which  $\sum\frac{1}{\lambda_n}<\infty$.
 We find the order of approximation of functions $f$ for which
 $\int_X|\lambda( {\chi})\hat{f}(\chi)|^2d\nu(\chi)<\infty$  \\
   Bibliography: 32 titles.
\end{abstract}
\noindent
 keywords: local fields, refinable equation, tight wavelet frames, p-adic numbers, trees.
\Large

\section*{Introduction}

The field $\mathbb Q_p$ of $p$-adic numbers is the completion of the field $\mathbb Q$ of rational numbers with respect of $p$-adic norm $|\cdot|_p$.
The notion of p-adic MRA was introduced and a general scheme
for its construction was described in \cite{ShSk}. The first orthogonal wavelet basis was con\-structed  by S. Kozyrev \cite{SK}. A. Yu. Khrennikov, V. M. Shelkovich, M.A.Skopina \cite{KhShSk} described a wide class of ortho\-gonal scaling functions generating  p-adic MRA.   An application p-adic wavelets  to pseudo-differential operators   are considered in \cite{KhSh2}.  A detailed exposition of p-adic analysis can be found in \cite{AKS}.
Non-compactly supported p-adic wavelets  were considered in \cite{E1}.

 S.Albeverio,  S.Evdokimov and M.Skopina \cite{AES} proved:
  1)there exists a unique p-adic MRA generated by an orthogonal scaling test function,
  2) there exists non\-orthogonal p-adic MRA
 3) there exists p-adic wavelet frames.     \\
  S.Evdokimov and M.Skopina \cite{ES} proved that any orthogonal $p$-adic  wavelet basis consisting of band-limited (periodic) functions is a modification of Haar basis.
In \cite{FLS}   an explicit description of all Vilenkin polynomials generating tight wavelet frames on Vilenkin group are given and an approximation error for functions from the Sobolev space with power weight function is obtained.

    In artcle \cite{LV1}, methods for constructing  $p$-adic step tight wavelet frames are indicated, but approximation estimates are missing.
       In this article we  obtain the $p$-adic  approximation estimates  by means of  step tight  wavelet frames, constructed in \cite{LV1} . As a corollary, we obtain an approximation estimate for functions from Sobolev spaces with logarithmic weight.  Note that in the classical case, a similar approximation estimate for functions from Sobolev spaces with power  weight, was obtained in  \cite{DHRS}.

        The problem of constructing tight wavelet frames  in zero-dimensional groups and the problem of approximation by means of these frames consists of two parts. Firstly, we must construct a scaling function, secondly, we must construct a tight
         frame from this scaling function, and, thirdly, we must obtain approximation estimates. The methods for solving the second and third problems are the same on any zero-dimensional group. But the methods for constructing scaling functions in different zero-dimensional groups are different and depend on the operation in the zero-dimensional group. Therefore, we consider all constructions and proofs in the zero-dimensional group. When constructing the scaling function, we assume that the addition operation in the zero-dimensional group satisfies the condition $pg_n=g_{n+1}$

\section{Notations and basic facts  }

Let $\mathbb Q_p$ be a field  of $p$-adic numbers \cite{AKS},  $\mathbb Q_p^+$ be an additive subgroup .
 Every $p$-adic number $x$ is the sum of the series $x=\sum_{k=m}^\infty x_k p^k$ in $p$-adic norm, where $x_k=0,1,...,p-1$. If $x_m\ne 0$, then $|x|_p=p^{-m}$.
 Subgroups $G_n=\{x=\sum_{k=n}^\infty x_k p^k\}$ form a decreasing sequence and define the base of the topology.
This means that the additive subgroup $\mathbb Q _p^+$ of the field $\mathbb Q _p$ is a locally compact zero-dimensional group
and  elements $g_k =p^k$ form a basis. It's obvious that $pg_k=g_{k+1}$.
In what follows, we will use only the additive group of the field $\mathbb Q _p$.

   Let us present the necessary information from the theory of zero-dimensional groups \cite{AVDR}. Let $(G,\dot + )$~be a~locally compact zero-dimensional Abelian
group with the topology generated by a~countable system of open
subgroups
\begin{equation}
\label{eq1.001}
\cdots\supset G_{-n}\supset\cdots\supset G_{-1}\supset G_0\supset
G_1\supset\cdots\supset G_n\supset\cdots
\end{equation}
where
$$
\bigcup_{n=-\infty}^{+\infty}G_n= G,\quad
 \quad \bigcap_{n=-\infty}^{+\infty}G_n=\{0\},
$$
 $p$ is an order of quotient groups $G_n/G_{n+1}$ for all $n\in\mathbb Z$.
 We will name such chain
as \it basic chain. \rm  In this case, a~base of the topology is
formed by all possible cosets~$G_n\dot + g$, $g\in G$.
The collection of such cosets, along with the empty set, forms a~semiring.
Given
a~coset $ G_n\dot+g$, we define a~measure~$\mu$ on it by
$\mu( G_n\dot+g)=\mu( G_n)= p^{-n}$. The measure~$\mu$ can be
extended onto the $\sigma$-algebra of measurable sets in the
standard way.
Given $n\in\mathbb Z$, consider an element $g_n\in G_n\setminus
G_{n+1}$ and fix~it. Then any $g\in G$ has a~unique representation
in the form
\begin{equation}
\label{eq1.01} g=\sum_{n=-\infty}^{+\infty}a_ng_n, \qquad
a_n=\overline{0,p-1}.
\end{equation}
The sum~\ref{eq1.01} contains finite number of terms with negative
subscripts, that~is,
\begin{equation}
\label{eq1.02} g=\sum_{n=m}^{+\infty}a_ng_n, \qquad
a_n=\overline{0,p-1}, \quad a_m\ne 0.
\end{equation}
We will name the system $(g_n)_{n\in \mathbb Z}$ as {\it a basic
system}.

By $X$ we denote  the collection of the characters of a~group $(G,\dot+ )$; it is
a~group with respect to multiplication, too. Also let
$G_n^\bot=\{\chi\in X:\forall\,x\in G_n\  , \chi(x)=1\}$ be the
annihilator of the group~$G_n$. Each annihilator~$ G_n^\bot$ is
a~group with respect to multiplication, and the subgroups~$
G_n^\bot$ form an~increa\-sing sequence
\begin{equation}
\label{eq1.03} \cdots\subset G_{-n}^\bot\subset\cdots\subset
G_0^\bot \subset G_1^\bot\subset\cdots\subset
G_n^\bot\subset\cdots
\end{equation}
with
$$
\bigcup_{n=-\infty}^{+\infty} G_n^\bot=X \quad {and} \quad
\bigcap_{n=-\infty}^{+\infty} G_n^\bot=\{1\},
$$
the quotient group $ G_{n+1}^\bot/ G_n^\bot$ having order~$p$.
The group of characters~$X$ is a zero-dimensional group with a
basic chain \ref{eq1.03}.  The family of
the cosets $ G_n^\bot\cdot\chi$, $\chi\in X$ form a base of topology.  Given
a~coset $ G_n^\bot\cdot\chi$, we define a~measure~$\nu$ on it by
$\nu( G_n^\bot\cdot\chi)=\nu( G_n^\bot)= p^n$. The measure~$\nu$ can be
extended onto the $\sigma$-algebra of measurable sets in the
standard way. One then forms the absolutely convergent integral
$\displaystyle\int_XF(\chi)\,d\nu(\chi)$ using this measure.

The value~$\chi(g)$ of the character~$\chi$ at an element $g\in G$
will be denoted by~$(\chi,g)$. The Fourier transform~$\widehat f$
of an~$f\in L_2( G)$  is defined~as follows
$$
\widehat f(\chi)=\int_{ G}f(x)\overline{(\chi,x)}\,d\mu(x)=
\lim_{n\to+\infty}\int_{ G_{-n}}f(x)\overline{(\chi,x)}\,d\mu(x),
$$
with the limit being in the norm of $L_2(X)$. For any~$f\in L_2(G)$,
the inversion formula is valid
$$
f(x)=\int_X\widehat f(\chi)(\chi,x)\,d\nu(\chi)
=\lim_{n\to+\infty}\int_{ G_n^\bot}\widehat
f(\chi)(\chi,x)\,d\nu(\chi);
$$
here the limit also signifies the convergence in the norm of~$L_2(
G)$. If $f,g\in L_2( G)$ then the Plancherel formula is valid \cite{AVDR}
$$
\int_{ G}f(x)\overline{g(x)}\,d\mu(x)= \int_X\widehat
f(\chi)\overline{\widehat g(\chi)}\,d\nu(\chi).
$$
\goodbreak

Provided with this topology, the group of characters~$X$ is
a~zero-dimensional locally compact group; there is, however,
a~dual situation: every element $x\in G$ is a~character of the
group~$X$, and~$ G_n$ is the annihilator of the group~$ G_n^\bot$.
We will denote the union of disjoint sets $E_j$  by $\bigsqcup E_j$.

 For any $n\in \mathbb Z$ we choose a character $r_n\in  G_{n+1}^{\bot}\backslash G_n^{\bot}$
 and fixed it. The collection of functions $(r_n)_{n\in \mathbb Z}$ is called a Rademacher system. Any character $\chi$ can be rewritten as a product
 $$\chi=\prod_{j=-m}^{+\infty} r_j^{\alpha_j},\ \alpha_j=\overline{0,p-1}.$$

  Let us denote
  $$
    H_0=\{h\in G: h=a_{-1}g_{-1}\dot+a_{-2}g_{-2}\dot+\dots \dot+ a_{-s}g_{-s}, s\in \mathbb
    N,\ a_j=\overline{0,p-1}\},
  $$
  $$
    H_0^{(s)}=\{h\in G: h=a_{-1}g_{-1}\dot+a_{-2}g_{-2}\dot+\dots \dot+
    a_{-s}g_{-s},\ a_j=\overline{0,p-1}
    \},s\in \mathbb N.
  $$
  In  $p$-adic group $\mathbb Q_p^+ $ the set
  $ H_0$    is usually denoted by $I_p$.
  Thus the set $H_0$ is an analog of the set $\mathbb N_0$.
  \begin{definition}
 We define the mapping ${\cal A}\colon G\to G$ by
 ${\cal A}x:=\sum_{n=-\infty}^{+\infty}a_ng_{n-1}$, where
 $x=\sum_{n=-\infty}^{+\infty}a_ng_n\in G$. As any element $x\in G$ can
 be uniquely expanded as~$x=\sum a_ng_n$, the mapping ${\cal A}\colon G\to
 G$ is one-to-one onto. The mapping~${\cal A}$ is called
 a  dilation operator if~${\cal A}(x\dot+ y)={\cal A}x\dot + {\cal A}y$ for all
 $x,y\in G$.
 \end{definition}

     By definition, put $(\chi {\cal A},x)=(\chi, {\cal
 A}x)$. It is also clear that
 ${\cal A} g_n= g_{n-1}, r_n{\cal A} = r_{n+1}$,\
 ${\cal A} G_n= G_{n-1}, G_n^\bot{\cal A}= G_{n+1}^\bot$.

 \begin{definition}[\cite{LS2}]
Let $M,N\in\mathbb N$.
We denote by  ${\mathfrak D}_M(G_{-N})$ the set of functions
 $f\in L_2(G)$ such that 1) ${\rm supp}\,f\subset G_{-N}$, and 2)
 $f$ is constant on cosets $G_M\dot+g$. The class ${\mathfrak
 D}_{-N}(G_{M}^\bot)$ is defined similarly.
 \end{definition}

 For a given function $\varphi\in  {\mathfrak D}_M(G_{-N})$ , we define  subspaces $V_n\subset L_2(G)$ generated by $\varphi$ as
  $$
 V_n=\overline{{\rm span}\{\varphi(\mathcal{A}^n\cdot \dot-h),h\in H_0\}}, \ n\in\mathbb Z.
 $$

 We say that sequence of subspaces $\{V_n\}$ forms a {\it multiresolution analysis} (MRA) for $L_2(G)$, if the folloving conditions are satisfied
 \begin{equation}\label{eq1.08}
 V_n\subset V_{n+1}, n\in \mathbb Z,
 \end{equation}
 \begin{equation}\label{eq1.09}
 \overline{\cup_n V_n}=L_2(G),
 \end{equation}
  \begin{equation}\label{eq1.1}
 \cap_n V_n=\{0\}
 \end{equation}

A function $\varphi \in L_2(G)$  is called refinable if
  \begin{equation}                                      \label{eq1.2}
   \varphi(x)=p\sum_{h\in H_0}\beta_h\varphi({\cal
   A}x\dot-h),
   \end{equation}
      for some sequence  $(\beta_h)\in \ell^2$. The equality   (\ref{eq1.2}) is called the refinement  equation. In Fourier domain, the equality   (\ref{eq1.2})  can be rewritten as
  $$
 \hat\varphi(\chi)=m_0(\chi)\hat\varphi(\chi{\cal A}^{-1}),
 $$
 where
 \begin{equation}                                           \label{eq1.3}
 m_0(\chi)=\sum_{h\in
 H_0}\beta_h\overline{(\chi{\cal A}^{-1},h)}
 \end{equation}
   is a mask of  (\ref{eq1.2}).

    \begin{lemma}[\cite{LVd2}]
    If a refinable function  $\varphi\in \mathfrak{D}_{G_{M}}(G_{-N}),\ M,N\in \mathbb N$, then its refinement equation is
      \begin{equation}                                      \label{eq1.4}
   \varphi(x)=p\sum_{h\in H_0^{(N+1)}}\beta_h\varphi({\cal
   A}x\dot-h),
   \end{equation}
   \end{lemma}

 If the shift system
 $(\varphi(x\dot-h))_{h\in H_0}$ form an  orthonormal basis in $V_0$, then MRA $(V_n)$ is called orthogonal.  Orthogonal MRA is used to construct orthogonal affine systems that form a basis of $L_2(\mathbb G)$.

  If the conditions  $\hat{\varphi} ( G_{-N}^\bot)=1$ and
  $|\hat\varphi(\chi)|\le 1$ are satisfied for    the refinable function $\varphi\in \mathfrak{D}_{G_{0}}(G_{-N})$, then  the function $\varphi$ generates an  MRA.  If $|\hat\varphi(\chi)|={\bf 1}_{G_0^\bot}(\chi)$, then this multiresolution  analysis is orthogonal \cite{LSF}.

The construction of tight frame systems starts with the construction of the system $\Psi=\{\psi_1,...,\psi_q\}\subset L_2(G)$. The objective  of MRA-construction  of tight  wavelet based frames is to find $\Psi=\{\psi_1,...,\psi_q\}\subset  V_1$ such that
$$
L(\Psi):=\{\psi_{\ell,n,h}=p^\frac{n}{2}\psi_\ell(\mathcal{A}^n\cdot \dot - h):1\le\ell \le q; n\in\mathbb{Z}, h\in H_0\}
$$
 forms a tight  frame for $L_2(G)$. The system $L(\Psi)\subset L_2(G)$ is called a {\it tight wavelet frame}
of $ L_2(G)$ if
$$
\|f\|^2_{L_2(G)}=\sum\limits_{g\in L(\Psi)}|\langle f,g\rangle|^2,
$$
holds for all $f\in L_2(G)$, where $\langle \cdot,\cdot\rangle$ is the inner product in $L_2(G)$. It is equivalent to
$$
f=\sum\limits_{g\in L(\Psi)}\langle f,g\rangle g,
$$
 for all $f\in L_2(G) $.
 Since $V_1$ is a $H_0$-invariant subspace generated by $\varphi(\mathcal{A}\cdot)$, finding $\Psi\subset V_1$ is the same as finding $\beta_h^{(\ell)}$ such that
 $$
 \psi_\ell (x)=p\sum\limits_{h\in H_0}\beta_h^{(\ell)}\varphi(\mathcal{A}x\dot- h).
 $$
 In  Fourier domain, this equality     can be rewritten as
  $$
 \hat\psi_\ell(\chi)=m_\ell(\chi)\hat\varphi(\chi{\cal A}^{-1}).
 $$

We will consider only step refinable function
$\varphi(x)\in \mathfrak{D}_{G_{M}}(G_{-N})$, that is equivalent to
$\hat\varphi(\chi)\in \mathfrak{D}_{G_{-N}^\bot}(G_M^\bot),\ M,N\in\mathbb N$.

If the shifts $(\varphi(x\dot-h))_{h\in H_0}$  are not orthogonal, then one can try to choose the functions $\psi^{(\ell)}(x)$  so that for any $ f\in L_2(G)$
   $$
   f(x)=\sum_{\ell =1}^r\sum_{n\in \mathbb Z}\sum_{h\in H_0}(f,\psi^{(\ell)} ({\cal A}^n \cdot\dot - h))\psi^{(\ell)} ({\cal A}^n x\dot - h).
   $$
    Such  system is called Parseval   wavelet frame or tight wavelet frame.

\section{Constructing refinable functions and tight wavelet frames}
In this section we will construct a  refinable
 function in the group $G$  which is the additive group of the field $\mathbb{Q}_p$. This suggests that
  $pg_n=g_{n+1}$.
 In this case $(r_n,g_m)=e^{\frac{2\pi i}{p^{n-m+1}}}$.

Let $M,N\in \mathbb N$.  We want to construct a refinable function
 $\varphi\in\mathfrak{ D}_{G_M}({G}_{-N})$, i.e.
$\hat{\varphi}\in\mathfrak{ D}_{G_{-N}^\bot}({G}_{M}^\bot)$.
Let us write the refinable equation in the Fourier domain

 \begin{equation}\label{eq2.1}
 \hat\varphi_(\chi)=\hat\varphi(\chi A^{-1})m_0(\chi).
  \end{equation}
 The mask
  \begin{equation} \label{eq2.2}
  m_0(\chi)=\sum_{h\in
  H_0^{(N+1)}}\beta_h\overline{(\chi,A^{-1}h)}
 \end{equation}
 is constant on cosets  ${G}_{-N}^\bot r_{-N}^{\alpha_{-N}}...r_{-N+s}^{\alpha_{-N+s}}$ .
 Let $m_0({G}_{-N}^\bot )=1$. Then
$$
\hat{\varphi}(\chi)=m_0(\chi)m_0(\chi\mathcal A^{-1})...m_0(\chi\mathcal A^{-N-M}).
$$
Denote by
 \begin{equation} \label{eq2.3}
 m_0({G}_{-N}^\bot r_{-N}^{\alpha_{-N}}r_{-N+1}^{\alpha_{-N+1}}...r_{0}^{\alpha_{0}}...r_{M}^{\alpha_{M}})=: \lambda_{ \alpha_{-N}\alpha_{-N+1}...\alpha_{0}...\alpha_{M}}=\lambda_m,
 \end{equation}
 where
 $$
 m= \alpha_{-N}+\alpha_{-N+1}p+...+\alpha_{0}p^N +...+\alpha_M p^{N+M} ,
 $$
  the values of the mask on  ${G}_{M+1}^\bot$. Since $(\chi \mathcal{A}^{-1},h)$ is constant on cosets
 $$
 {G}_{-N}^\bot r_{-N}^{\alpha_{-N}}r_{-N+1}^{\alpha_{-N+1}}...r_{M}^{\alpha_{M}}
$$
we have
$$
(\chi\mathcal{A}^{-1}, h)=
\prod\limits_{\nu=0}^{-N} \prod\limits_{k=-N}^M e^{\frac{2\pi i}{p^{k-\nu+1}}\alpha_k a_{\nu}}.
$$

Therefore
$$
m_0(\chi)
=\sum\limits_{a_{-1},...,a_{-N-1}}\beta_{a_{-1},...,a_{-N-1}}
e^{-\frac{2\pi i}{p}\sum\limits_{\nu=0}^{-N}a_{\nu-1}p^{\nu+N}\sum\limits_{k=-N}^M \alpha_k p^{-k
-N} }.
$$
Denote

$$
q_m=e^{-\frac{2\pi i}{p}}\sum\limits_{k=-N}^{M}\alpha_k p^{-N-k}, \ n=a_{-N-1}+a_{-N}p+...+a_{-1}p^N.
$$
and write the equality  (\ref{eq2.2})  in the form

\begin{equation}\label{eq2.4}
 \left(\begin{array}{rrrr}
q_{0}^0&q_{0}^1&...&q_{0}^{p^{N+1}-1}\\
q_{1}^0&q_{1}^1&...&q_{1}^{p^{N+1}-1}\\
...&...&...&...\\
q_{p^{N+1}-1}^0&q_{p^{N+1}-1}^1&...&q_{p^{N+1}-1}^{p^{N+1}-1}\\
...&...&...&...\\
q_{p^{M+N+1}-1}^0&q_{p^{M+N+1}-1}^1&...&q_{p^{M+N+1}-1}^{p^{N+1}-1}\\

\end{array}
\right)
\left(\begin{array}{c}
\beta_0\\
\beta_1\\
\beta_2\\
\vdots\\
\beta_{p^{N+1}-1}
\end{array}
\right)=
\left(\begin{array}{c}
\lambda_0\\
\lambda_1\\
\lambda_2\\
\vdots\\
\lambda_{p^{M+N+1}-1}\\
\end{array}
\right)
 \end{equation}

We need to find $\lambda_m$ and $\beta_n$ so that
$$
\hat{\varphi}(\chi)=m_0(\chi)m_0(\chi\mathcal A^{-1})...m_0(\chi\mathcal A^{-N-M})=0
$$
on the set ${G}_{M+1}^\bot\setminus {G}_{M}^\bot$.
To find $\lambda_m$  and $\beta_n$ we construct a rooted $p$-{\it adic} mask tree $T=T(m_0)$ in the following way.
For any $m\in \mathbb N: p^{M+N}\le m\le p^{M+N+1}-1$ we construct the path
$$
\lambda_m\rightarrow \,\lambda_{m \,{\rm div}\, p} \rightarrow \lambda_{m \,{\rm div}\, p^2}\dots \rightarrow \lambda_{m \,{\rm div}\, p^{M+N}}\rightarrow \lambda_0=1.
$$

 from the leaf $\lambda_m$ to the root $\lambda_0$ of the  tree $T(m_0)$, where $a \,{\rm div}\, b$ stands for the integer division of $a$ by $b$. It is clear that $H=M+N$ is a high of this $p$-{\it adic} tree $T$. (See Figure 2 for graph of $T$). In this tree, the numbers  $\lambda_{j_s}: p^{s-1}\le j_s\le p^{s}-1, s\ge 0$ form the s-th level. The set of all products $\lambda_{m}\lambda_{m \,{\rm div} \, p}...\lambda_0$ coincides with the set of all values of the function $\hat{\varphi}(\chi)$ on
 the set ${G}_{M+1}^\bot\setminus {G}_{M}^\bot$. On each path
$$
\lambda_m\rightarrow \,\lambda_{m \,{\rm div}\, p} \rightarrow \lambda_{m \,{\rm div}\, p^2}\dots \rightarrow \lambda_{m \,{\rm div}\, p^{M+N}}
$$
 we select one node and place zero there.  Denote the set of chosen zeros by  $\Lambda_0 (T)$.

\unitlength=0.80mm
  \begin{picture}(240,80)
 \small

  \put(150,48){\line (-1,1){16}}
     \put(151,48){\line (0,1){16}}
     \put(152,48){\line (1,1){16}}
     \put(125,66){$\lambda_{p^{N+M+1}-p}$}
      \put(160,66){$\lambda_{p^{N+M+1}-1}$}

     \put(58,48){\line (-1,1){16}}
     \put(59,48){\line (0,1){16}}
     \put(60,48){\line (1,1){16}}
     \put(38,66){$\lambda_{p^{N+M}}$}
      \put(70,66){$\lambda_{p^{N+M}+p-1}$}
      \multiput(51,64)(3,0){6}{$\cdot$}

    \multiput(60,46)(3,0){31}{$\cdot$}
   \put(130,25){\line (-1,1){16}}
     \put(131,25){\line (0,1){16}}
     \put(132,25){\line (1,1){16}}
     \put(110,43){$\lambda_{p^{N}-p}$}
      \put(146,43){$\lambda_{p^{N}-1}$}
      \multiput(124,41)(3,0){6}{$\cdot$}
      \multiput(74,41)(3,0){4}{$\cdot$}

  \put(124,19){$\lambda_{p-1}$}
  \put(81,19){$\lambda_{1}$}
   \put(105,19){$\lambda_{\nu}$}
   \put(80,25){\line (-1,1){16}}
     \put(81,25){\line (0,1){16}}
     \put(82,25){\line (1,1){16}}
     \put(60,43){$\lambda_{p^{N-1}}$}
      \put(86,43){$\lambda_{p^{N-1}+p-1}$}

  \multiput(81,22)(3,0){17}{$\cdot$}

 \put(100,-2){$\lambda_{0}=1$}
  \put(100,2){\line (-1,1){16}}
  \put(106,2){\line (0,1){16}}
  \put(110,2){\line (1,1){16}}
  \put(-5,66){$\bf  {G}_{M+1}^\bot \setminus {G}_{M}^\bot $:}
  \put(-5,43){$\bf {G}_{0}^\bot \setminus {G}_{-1}^\bot $:}
  \put(-5,19){$\bf {G}_{-N+1}^\bot \setminus {G}_{-N}^\bot $:}
  \put(-5,-2){$\bf {G}_{-N}^\bot  $:}

    \end{picture}\\

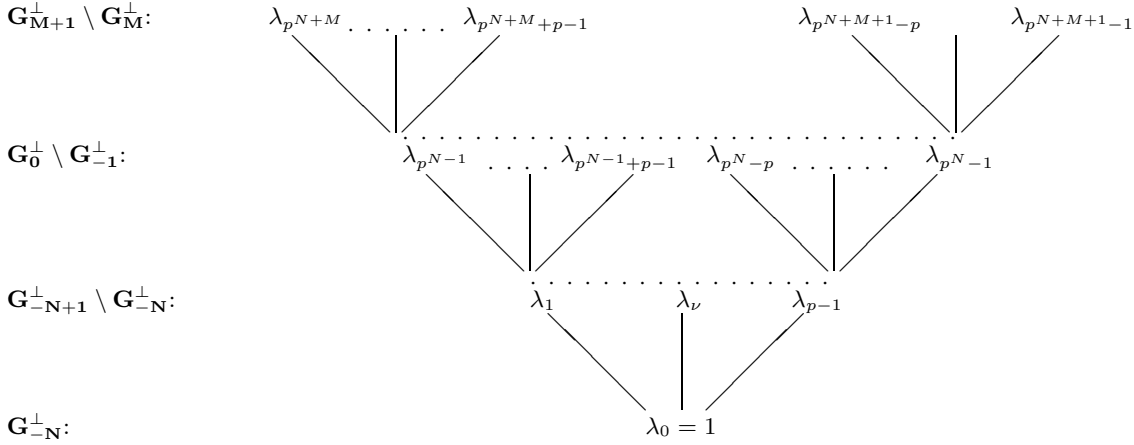
\hskip4cm Figure 2.  The graph of the tree $T=T(m_0)$.

\begin{theorem}[\cite{LV1}]\label{th2.01}
 1) Let $\sharp\Lambda_0(T)) = p^{N+1}-1$. Then  the corresponding  values $\lambda_\nu\in \Lambda_0(T)$ determine the mask $m_0$ of some refinable  function. If $\lambda_\nu\notin \Lambda_0(T) $  then $\lambda_\nu\neq 0$.\\
 2) Let $\sharp\Lambda_0(T)< p^{N+1}-1$ . Then  the corresponding  values $\lambda_\nu\in \Lambda_0(T)$ determine the mask $m_0$ of some refinable  function.\\
 3) If $\sharp\Lambda_0(T)\ge p^{N+1}$ then  the corresponding  values $\lambda_\nu\in \Lambda_0(T)$ do not define a mask of the refinable  function.
 \end{theorem}

 {\bf Corollary.}\cite{LV1}Using this method of constructing a mask, we can obtain all masks that generate $\hat{\varphi}\in \mathfrak{D}_{G_{-N}^\bot}({G}_{M}^\bot))$ and hence all refinable  functions $\hat{\varphi}\in \mathfrak{D}_{G_{-N}^\bot}(G_{M}^\bot)$.\\

For the function $\varphi\in L_2({G})$ we will use the standard notation
$$
\varphi_{n,h}=p^{\frac{n}{2}}\varphi(\mathcal{A}^n\cdot \dot-h), \quad h\in H_0, n\in \mathbb Z.
$$
Let the mask  $m_0({\chi})$ built on a tree $T$ and $\varphi\in \mathfrak D_{G_{M}}({G}_{-N})$ be the corresponding refinable  function.
  We want to find masks $m_j, j=\overline{1,q}$ and corresponding functions $\psi^{(j)}$ that generate a tight wavelet frame. We will use lemma  \ref{Lm3.1} and theorem 2.2 which are valid for any zero-dimensional group $G$.

\begin{lemma}[\cite{LV1}]\label{Lm3.1}
Let the mask $m_j \ (j=1,...,q)$ satisfy the following conditions:
1)$
\hat{\psi}^{(j)}(\chi)=\hat{\varphi}(\chi\mathcal{A}^{-1})m_j(\chi)=
\xi_{\gamma_{-s},...,\gamma_{-1},\gamma_{0},...,\gamma_{u}}\in \mathbb C $   on the coset
 $$   G^\bot_{-s}r_{-s}^{\gamma_{-s}}
... r_{-1}^{\gamma_{-1}}r_{0}^{\gamma_{0}}...r_{u}^{\gamma_{u}}, \quad s=s(j), s\le N
$$
 and\\
 2) $m_j(\chi)=0$ outside the coset $   G^\bot_{-s}r_{-s}^{\gamma_{-s}}
... r_{-1}^{\gamma_{-1}}r_{0}^{\gamma_{0}}...r_{u}^{\gamma_{u}}$, \\
3)$ |\xi_{\gamma_{-s},...,\gamma_{-1},\gamma_{0},...,\gamma_{u}}|=1$.
Then
$$
\sum_{h\in H_0} |c_{n,h}^{(j)}(f)|^2=
\sum_{h\in H_0} |(\psi^{(j)}_{n,h},f)|^2=
\int_{G_{n-s}^\bot r_{n-s}^{\gamma_{-s}}...r_{n-1}^{\gamma_{-1}}r_{n-0}^{\gamma_{0}}...r_{n+u}^{\gamma_{u}}}|\hat{f}(\chi) |^2d\nu(\chi)
$$
\end{lemma}

\begin{theorem}[\cite{LV1}]\label{Th2.2}
  Let $\varphi \in \mathfrak{D}_{G_{M}}(G_{-N})$ be a refinable function with a mask $m_0$. Define masks  $m_j: j=1,2,...,q$ so that  \\
1)$\hat\varphi (\chi\mathcal{A}^{-1})m_j(\chi)=
{\bf 1}_{E_j}(\chi)$, where  $E_j=G^\bot_{-s(j)}r_{-s(j)}^{\alpha_{-s(j)}}r_{-s(j)+1}^{\alpha_{-s(j)+1}}...r_{0}^{\alpha_{0}}...r_{M}^{\alpha_{M}}$ are disjoint cosets and  $E_j \mathcal{A}^t$ are disjoint also ,\\
2)there  are integers   $t(j)\ge 0$, such  that
$$
\bigsqcup_jE_j\mathcal{A}^{t(j)}=G_{M+1}^\bot\setminus G_M^\bot.
$$
 Then functions   $
\psi^{(1)},\psi^{(2)},...,\psi^{(q)}$ generate tight wavelet frame .
\end{theorem}
\section{Approximation order }
In this section we study approximation properties of constructing frames on $p$-adic group $G$.
Let the mask $m_0$ be constructed  in according with the tree $T$, $\varphi \in \mathfrak{D}_{G_{M}}(G_{-N})$ be a refinable function with a mask $m_0$. Now we construct masks $m_j$ in accordance with Theorem 2.2. By Theorem 2.2, the functions $\{\psi_{n,h}^{(j)}\}$ form a tight frame.
Since $\{\psi_{n,h}^{(j)}\}$ form  a tight frame, then
$$\lim\limits_{\tilde{N}\to +\infty}\|f-  \sum_{n=-\infty}^{\tilde{N}}\sum_{j=1}^q\sum_{h\in H_0}\langle f,\psi_{n,h}^{(j)}\rangle \psi_{n,h}^{(j)}\|_2=0.
$$
We will find the order of this approximation.
\begin{theorem}Let the functions $\psi_1,\psi_2,...,\psi_q$ and numbers $t(j)$ be constructed as in the Theorem \ref{Th2.2}. Let us denote $l=\max t(j)$. Then for $\tilde{N}>N$  the inequalities
$$\|f-  \sum_{n=-\infty}^{\tilde{N}}\sum_{j=1}^q\sum_{h\in H_0}\langle f,\psi_{n,h}^{(j)}\rangle \psi_{n,h}^{(j)}\|_2
\le (N+1)p^{\frac{M-1}{2}}\sum_{n=\tilde{N}+1}^\infty \left(\int_{G_{n-l+1}^\bot\setminus G^\bot_{n-l}}|\hat{f}(\chi)|^2d\nu(\chi)\right)^\frac12 $$
are satisfied.
\end{theorem}
{\bf Proof.} By the property of tight frames  we have the inequality \cite{FLS}
$$
R_{\tilde{N}}\stackrel{\text{def}}=\|f-  \sum_{n=-\infty}^{\tilde{N}} \sum_{j=1}^q\sum_{h\in H_0}\langle f,\psi_{n,h}^{(j)}\rangle \psi_{n,h}^{(j)}\|_2=
\|\sum_{n=\tilde{N}+1}^\infty \sum_{j=1}^q \sum_{h\in H_0}\langle f,\psi_{n,h}^{(j)}\rangle \psi_{n,h}^{(j)}\|_2\le .
$$
$$
 \le\sum_{n=\tilde{N}+1}^\infty \|\sum_{j=1}^q\sum_{h\in H_0}\langle f,\psi_{n,h}^{(j)}\rangle \psi_{n,h}^{(j)}\|_2.
$$
Let us denote $c^{(j)}_{n,h}:=\langle f,\psi_{n,h}^{(j)}\rangle$. Using the equality $\hat{\psi}^{(j)}_{n,h}=\hat{\psi}^{(j)}(\chi\mathcal{A}^{-n})\overline{(\chi \mathcal{A}^{-n},h)}$ we have
$$
\|\sum_{j=1}^q\sum_{h\in H_0}\langle f,\psi_{n,h}^{(j)}\rangle \psi_{n,h}^{(j)}\|_2^2=\|\sum_{j=1}^q\sum_{h\in H_0} c^{(j)}_{n,h} \hat{\psi}_{n,h}^{(j)}\|_2^2=
$$
$$
=\int_X|\sum_{j=1}^q\sum_{h\in H_0} c^{(j)}_{n,h} \hat{\psi}_{n,h}^{(j)}(\chi)|^2d\nu(\chi)=
$$

$$=
\int_X|\sum_{j=1}^q\sum_{h\in H_0} c^{(j)}_{n,h} p^{-\frac{n}{2}}\hat{\psi}^{(j)}(\chi\mathcal{A}^{-n})
(\overline{\chi\mathcal{A}^{-n},h)}|^2d\nu(\chi)=
$$

$$
=\frac{1}{p^n}\int_X|\sum_{j=1}^q\sum_{h\in H_0} c^{(j)}_{n,h} \hat{\psi}^{(j)}(\chi\mathcal{A}^{-n})
(\overline{\chi\mathcal{A}^{-n},h)}|^2d\nu(\chi)=
$$

$$
=\int_X|\sum_{j=1}^q\sum_{h\in H_0} c^{(j)}_{n,h} \hat{\psi}^{(j)}(\chi)
(\overline{\chi,h)}|^2d\nu(\chi)=
$$
( we used the fact that $\hat{\psi}^{(j)}={\bf 1}_{E_j}$ and  sets  $E_j$ are disjoint)
$$
=\int_X\biggl|\sum_{j=1}^q\hat{\psi}^{(j)}(\chi)\biggr|^2 \biggl|\sum_{h\in H_0} c^{(j)}_{n,h}
(\overline{\chi,h)}\biggr|^2d\nu(\chi)
=\sum_{j=1}\int_{E_j} \biggl|\sum_{h\in H_0} c^{(j)}_{n,h}
(\overline{\chi,h)}\biggr|^2d\nu(\chi).
$$
If ${E_j}\subset G_0^\bot$, then

$$
\int_{E_j} \biggl|\sum_{h\in H_0} c^{(j)}_{n,h}
(\overline{\chi,h)}\biggr|^2d\nu(\chi)\le \int_{G_0^\bot} \biggl|\sum_{h\in H_0} c^{(j)}_{n,h}
(\overline{\chi,h)}\biggr|^2d\nu(\chi)=\sum_{h\in H_0} |c^{(j)}_{n,h}|^2.
$$
If ${E_j}\subset G_1^\bot\setminus G_0^\bot$, then
$$
\int_{E_j} \biggl|\sum_{h\in H_0} c^{(j)}_{n,h}
(\overline{\chi,h)}\biggr|^2d\nu(\chi)\le \int_{G_0^\bot r_0^\alpha  } \biggl|\sum_{h\in H_0} c^{(j)}_{n,h}
(\overline{\chi,h)}\biggr|^2d\nu(\chi)=\sum_{h\in H_0} |c^{(j)}_{n,h}|^2.
$$

If ${E_j}\subset G_m^\bot\setminus G_{m-1}^\bot, (1<m\le M)$, then
$$
\int_{E_j} \biggl|\sum_{h\in H_0} c^{(j)}_{n,h}
(\overline{\chi,h)}\biggr|^2d\nu(\chi)\le \int_{G_{m-1}^\bot r_{m-1}^\alpha  } \biggl|\sum_{h\in H_0} c^{(j)}_{n,h}
(\overline{\chi,h)}\biggr|^2d\nu(\chi)=
$$
$$
=\sum_{\alpha_0,...,\alpha_{m-2}}\int_{G_0^bot r_0^{\alpha_0}...r_{m-2}^{\alpha_{m-2}}r_{m-1}^{\alpha_{m-1}}}
\biggl|\sum_{h\in H_0} c^{(j)}_{n,h}
(\overline{\chi,h)}\biggr|^2d\nu(\chi)
$$
$$
=\sum_{\alpha_0,...,\alpha_{m-2}}\sum_{h\in H_0} |c^{(j)}_{n,h}|^2=
p^{m-1}\sum_{h\in H_0} |c^{(j)}_{n,h}|^2\le p^{M-1}\sum_{h\in H_0} |c^{(j)}_{n,h}|^2.
$$
So
$$
\|\sum_{j=1}^q\sum_{h\in H_0}\langle f,\psi_{n,h}^{(j)}\rangle \psi_{n,h}^{(j)}\|_2^2 \le
p^{M-1}\sum_{j=1}^q\sum_{h\in H_0} |c^{(j)}_{n,h}|^2.
$$
By lemma \ref{Lm3.1}
$$
\sum_{h\in H_0} |c^{(j)}_{n,h}|^2=\sum_{h\in H_0}|\langle\hat{\psi}^{(j)}_{n,h},f\rangle|^2=
\int_{E_j\mathcal{A}^n}|\hat{f}(\chi)|^2d\nu(\chi).
$$
Therefore
\begin{equation}\label{Eq4.3}
\|\sum_{j=1}^q\sum_{h\in H_0}\langle f,\psi_{n,h}^{(j)}\rangle \psi_{n,h}^{(j)}\|_2^2 \le
p^{M-1}\sum_{j=1}^q\int_{E_j\mathcal{A}^n}|\hat{f}(\chi)|^2d\nu(\chi).
\end{equation}
Let   $l=t(j)+\alpha(j)$. It's obvious that $t(j)\le N,\alpha(j)\le N$.

If $n=l$ then

$$
\sum_{j=1}^q\int_{E_j\mathcal{A}^l}|\hat{f}(\chi)|^2d\nu=
\sum_{j=1}^q\int_{E_j\mathcal{A}^{t(j)}\mathcal{A}^{\alpha(j)}}|\hat{f}(\chi)|^2d\nu\le
$$
$$
\le\sum_{j=1}^q\left(\int_{E_j\mathcal{A}^{t(j)}}|\hat{f}(\chi)|^2d\nu+
\int_{E_j\mathcal{A}^{t(j)}\mathcal{A}}|\hat{f}(\chi)|^2d\nu+...
+\int_{E_j\mathcal{A}^{t(j)}\mathcal{A}^N}|\hat{f}(\chi)|^2d\nu
\right)=
$$
$$
=\int_{G_1^\bot\setminus G_0^\bot}|\hat{f}(\chi)|^2d\nu+\int_{G_2^\bot\setminus G_1^\bot}|\hat{f}(\chi)|^2d\nu+...+\int_{G_{N+1}^\bot\setminus G_N^\bot}|\hat{f}(\chi)|^2d\nu.
$$
If $n>l$ then we have
$$
\sum_{j=1}^q\int\limits_{E_j\mathcal{A}^n}|\hat{f}(\chi)|^2d\nu\le \int_{G_{n-l+1}^\bot\setminus G_{n-l}^\bot}|\hat{f}(\chi)|^2d\nu+...+
\int\limits_{G_{n-l+N+1}^\bot\setminus G_{n-l+N}^\bot}|\hat{f}(\chi)|^2d\nu.
$$
Now we can write  inequality  (\ref{Eq4.3}) as
$$
\|\sum_{j=1}^q\sum_{h\in H_0}\langle f,\psi_{n,h}^{(j)}\rangle \psi_{n,h}^{(j)}\|_2\le
p^{\frac{M-1}{2}}\sum_{k=0}^N \left(\int_{G_{n-l+1+k}^\bot\setminus G_{n-l+k}^\bot}|\hat{f}(\chi)|^2d\nu(\chi)\right)^\frac12.
$$
Summing up these inequalities, we obtain for $\tilde{N}>N$
$$
\sum_{n=\tilde{N}+1}^\infty\|\sum_{j=1}^q\sum_{h\in H_0}\langle f,\psi_{n,h}^{(j)}\rangle \psi_{n,h}^{(j)}\|_2\le
p^{\frac{M-1}{2}}\sum_{n=\tilde{N}+1}^\infty\sum_{k=0}^N \left(\int_{G_{n-l+1+k}^\bot\setminus G_{n-l+k}^\bot}|\hat{f}(\chi)|^2d\nu(\chi)\right)^\frac12.
$$
Each term on the right side of this inequality is present at most $N+1$ times. Therefore
\begin{equation}\label{Eq4.4}
\sum_{n=\tilde{N}+1}^\infty \|\sum_{j=1}^q\sum_{h\in H_0}\langle f,\psi_{n,h}^{(j)}\rangle \psi_{n,h}^{(j)}\|_2\le
(N+1)p^{\frac{M-1}{2}}\sum_{n=\tilde{N}+1}^\infty \left(\int_{G_{n-l+1}^\bot\setminus G_{n-l}^\bot}|\hat{f}(\chi)|^2d\nu\right)^\frac12.
\end{equation}
This completes the proof.\\
{\bf Remark.} From Theorem 3.1 we can obtain an approximation estimate for functions from  Sobolev spaces.  Let's choose an increasing sequence $\gamma_n\uparrow +\infty$ so that $\sum_{n=0}^{+\infty} \frac{1}{\gamma_n}<\infty$, and $ \gamma_n=1$ for   $n\in-\mathbb{N}$. Define the function $\gamma(\chi)=\gamma_n=\gamma(|\chi|_p)$ for $\chi\in G_{n}^\bot\setminus G_{n-1}^\bot$.
Let's transform the right side in the inequality (\ref{Eq4.4}).
$$
\sum_{n=\tilde{N}+1}^\infty \left(\int_{G_{n-l+1}^\bot\setminus G_{n-l}^\bot}|\hat{f}(\chi)|^2d\nu\right)^\frac12=
\frac{1}{p^{\frac{l}{2}}}\sum_{n=\tilde{N}+1}^\infty \left(\int_{G_{n+1}^\bot\setminus G_{n}^\bot}|\hat{f}(\chi \mathcal{A}^{-l})|^2d\nu\right)^\frac12=
$$
$$
=\frac{1}{p^{\frac{l}{2}}}\sum_{n=\tilde{N}+1}^\infty \frac{1}{\gamma_{n+1}}\left(\int_{G_{n+1}^\bot\setminus G_{n}^\bot}\gamma_{n+1}^2|\hat{f}(\chi \mathcal{A}^{-l})|^2d\nu\right)^\frac12=
$$
$$=\frac{1}{p^{\frac{l}{2}}}\sum_{n=\tilde{N}+1}^\infty \frac{1}{\gamma_{n+1}}\left(\int_{G_{n+1}^\bot\setminus G_{n}^\bot}\gamma^2(\chi)|\hat{f}(\chi \mathcal{A}^{-l})|^2d\nu\right)^\frac12\le
$$
$$
\le \frac{1}{p^{\frac{l}{2}}}\sum_{n=\tilde{N}+1}^\infty\frac{1}{\gamma_{n+1}} \left(\int_X \gamma^2(\chi)|\hat{f}(\chi \mathcal{A}^{-l})|^2d\nu\right)^\frac12.
$$
Therefore
$$
R_{\tilde{N}}\le \frac{(N+1)p^{\frac{M-1}{2}}}{p^{\frac{l}{2}}}\sum_{n=\tilde{N}+1}^\infty\frac{1}{\gamma_{n+1}} \left(\int_X \gamma^2(\chi )|\hat{f}(\chi \mathcal{A}^{-l})|^2d\nu\right)^\frac12=
$$
$$
=(N+1)p^{\frac{M-1}{2}} \left(\int_X \gamma^2(\chi \mathcal{A}^{l})|\hat{f}(\chi )|^2d\nu\right)^\frac12 \sum_{n=\tilde{N}+1}^\infty\frac{1}{\gamma_{n+1}}.
$$

\noindent
 If $G$ is  $p$-adic group and $\gamma_k=p^{km}\ (k\ge 0)$ then we obtain analog of Theorem 14 from \cite{FLS}
 \begin{equation}\label{Eq4.5}
R_{\tilde{N}}\le \frac{(N+1)p^{\frac{M-1}{2}}}{(p^m-1)p^{m{\tilde{N}+1}}} \left(\int_X (1+|\chi|_p^{m+l})^2|\hat{f}(\chi )|^2d\nu\right)^\frac12.
\end{equation}
Indeed, if $\chi\in G_k^\bot\setminus G_{k-1}^\bot$ then
$$
\gamma(\chi)=\gamma_k=
\left\{\begin{array}{lr}
p^{km}=|\chi|_p^m,& if \ k\ge 0,\\
1,&if \  k<0.\\
\end{array}\right.
$$
Therefore $\gamma(\chi)=\max(1,|\chi|_p^m)\le (1+|\chi|_p^m)$ and
$$
\int_X\gamma^2(\chi)|\hat{f}(\chi\mathcal{A}^{-l})|^2d\nu\le \int_X(1+|\chi|_p^m)^2|\hat{f}(\chi\mathcal{A}^{-l})|^2d\nu=
$$
$$
=p^l\int_X(1+|\chi|_p^{m+l})^2|\hat{f}(\chi)|^2d\nu.
$$
If we take $\gamma_k=(k+1)^{1+\varepsilon/2},\ (\varepsilon>0, k\ge 0) $ then
$$
\gamma^2(\chi)\le (1+\log_p^+|\chi|_p)^{2+\varepsilon},
$$
 and we get the inequality
  $$
R_{\tilde{N}}\le \frac{2(N+1)}{\varepsilon (1+\tilde{N})^{\varepsilon/2}} \left(\int_X \bigl(1+l+\log_p^+|\chi|_p\bigr)^{2+\varepsilon}|\hat{f}(\chi )|^2d\nu\right)^\frac12.
$$
In these inequalities $|\chi |_p= p^n$ for $\chi\in G_n^\bot\setminus G_{n-1}^\bot$, and
$$
\log_p^+|\chi|_p=
\left\{\begin{array}{lr}
\log_p|\chi|_p,& if\ |\chi|_p>1,\\
1,&if \ |\chi|_p \le 1.\\
\end{array}\right.
$$

\noindent
 {\bf Funding}\\This work was supported by the Russian Science Foundation
 No  22-21-00037,   https://rscf.ru/project/22-21-00037.
\vskip1cm
  \noindent
  {\bf ORCID}  http://orcid.org/0000-0003-3038-2698

 \end{document}